\documentclass[12pt]{article}
\usepackage{mathrsfs}
\usepackage{t1enc}
\usepackage[latin1]{inputenc}
\usepackage[english]{babel}
\usepackage{amsmath,amsthm}
\usepackage{amsfonts}
\usepackage{latexsym}
\usepackage{indentfirst}
\usepackage{listings}
\usepackage{hyperref}
\usepackage[overload]{empheq}
\usepackage{cleveref}
\usepackage[natural]{xcolor}

\usepackage{graphicx, epsfig, subfigure}
\usepackage{epstopdf}
\usepackage{enumerate} 
\usepackage[subnum]{cases}
\usepackage{geometry}
\usepackage{caption}

\theoremstyle{plain}

\newtheorem{theorem}{Theorem}
\numberwithin{theorem}{section}
\newtheorem{lemma}{Lemma}
\numberwithin{lemma}{section}

\numberwithin{corollary}{section}

\newtheorem{conjecture}{Conjecture}

\numberwithin{conjecture}{section}
\newtheorem*{conjecture*}{Conjecture}

\newtheorem{question}{Question}

\theoremstyle{definition}

\theoremstyle{remark}

\theoremstyle{plain}

\DeclareMathOperator{\md}{mod}

\theoremstyle{plain}

\theoremstyle{plain}

\usepackage{lineno}
\newcommand{\pf}{\noindent\textbf{Proof}.\quad}

\title{Antimagic orientation of biregular bipartite graphs\thanks{This work was supported by the National Natural Science Foundation of China (11371355,
11471193, 11271006, 11631014), the Foundation for Distinguished Young Scholars of Shandong Province
(JQ201501).}}
\author{Songling Shan$^a$\thanks{Email: songling.shan@vanderbilt.edu},\ \ Xiaowei Yu$^{a,b}$\thanks{Corresponding author. Email: xwyu2013@163.com, yux6@vanderbilt.edu}
\unskip\\[.5em]
{\small  $^a$ Department of Mathematics, Vanderbilt University,}\\
{\small  Nashville, TN 37240, U.S.A.}\\
{\small  $^b$ School of Mathematics, Shandong University,}\\
{\small  Jinan, Shandong 250100, P. R. China}\\
}

\date{}
\begin{document}
\baselineskip 0.65cm

\maketitle

\begin{abstract}
An antimagic labeling
of a directed graph $D$ with $n$ vertices and $m$ arcs is a bijection from the set of arcs of
$D$ to the integers $\{1, \cdots,  m\}$ such that all $n$ oriented vertex sums are pairwise distinct,
where an oriented vertex sum is the sum of labels of all arcs entering that vertex minus the
sum of labels of all arcs leaving it. An undirected graph $G$ is said to
have an  antimagic orientation if  $G$ has  an orientation which  admits an antimagic
labeling. Hefetz, M{\"{u}}tze, and Schwartz  conjectured that
every connected  undirected graph admits an antimagic orientation.
In this paper, we support this conjecture by proving that every  biregular bipartite graph admits an antimagic orientation.
\end{abstract}

\bigskip


\section{Introduction}
Unless otherwise stated explicitly, all graphs considered are simple and finite.
 A {\it labeling\/} of a graph $G$  with $m$ edges is a bijection from $E(G)$ to a set $S$ of $m$ integers,
 and the {\it vertex sum\/} at a vertex $v\in V(G)$ is the sum of labels on the edges incident to $v$.
 If there are two vertices have same vertex sums in $G$, then we call them \emph{conflict}.
A labeling of $E(G)$ with  no conflicting vertex is called a {\it vertex distinguishable labeling\/}.
A  labeling  is {\it antimagic\/} if it is vertex distinguishable and $S=\{1,2,\cdots, m\}$.
A graph is {\it antimagic\/} if it has an antimagic labeling.

Hartsfield and Ringel~\cite{MR1282717} introduced antimagic labelings in 1990 and conjectured
that every connected graph other than $K_2$ is antimagic.
There have been significant progresses toward this conjecture. Let $G$ be a graph with $n$ vertices
other than $K_2$.
In 2004, Alon, Kaplan, Lev, Roditty, and Yuster~\cite{MR2096791} showed that
there exists a constant $c$ such that if $G$ has minimum degree
at least $ c \cdot log n$,
then $G$ is antimagic.  They also proved that $G$ is antimagic when the
maximum degree of $G$ is at least $n-2$, and they
proved that all complete multipartite graphs (other than $K_2$) are antimagic.
The latter result of Alon et al. was improved by
Yilma~\cite{MR3021347} in 2013.

Apart from the above results on dense graphs,
the antimagic labeling conjecture has been also verified for
regular graphs.
Started with Cranston~\cite{MR2478227} showing that
every bipartite regular graph is antimagic,
regular graphs of odd degree~\cite{MR3372337},
and finally all  regular graphs~\cite{MR3515572} were  shown
to be antimatic sequentially.
For more results on the antimagic labeling conjecture
for other classes of graphs, see~\cite{MR3527991, MR2174213, MR2682515,MR2510327}.

Hefetz, M{\"{u}}tze, and Schwartz~\cite{MR2674494} introduced the
variation of antimagic labelings, i.e., antimagic labelings on
directed graphs. An {\it antimagic\/} labeling of a directed graph with  $m$ arcs is
a bijection from the set of arcs  to the integers $\{1,...,m\}$ such that
any two oriented vertex sums are pairwise distinct, where an {\it oriented
vertex sum\/} is the sum of labels of all arcs entering that vertex minus
the sum of labels of all arcs leaving it.
A digraph is called {\it antimagic\/} if it admits an antimagic labeling. For an undirected graph $G$, if it has an orientation such that
the orientation is antimagic, then we say $G$ admits an {\it antimagic orientation\/}.
Hefetz et al. in the same paper posted the following problems.

\begin{question}[\cite{MR2674494}]\label{question1}
 Is every connected directed graph with at least 4 vertices antimagic?
\end{question}

\begin{conjecture}[\cite{MR2674494}]\label{antimagic-orientation}
Every connected  graph admits an antimagic orientation.
\end{conjecture}

Parallel to the results the on antimagic labelling conjecture,
Hefetz, M{\"{u}}tze, and Schwartz~\cite{MR2674494} showed
that every orientation of a dense graph is antimagic
and almost all regular graphs
have an antimagic orientation. Particulary,
they showed that every orientation of stars\,(other than $K_{1,2}$), wheels,
and complete graphs\,(other than $K_3$)
is antimagic. Observe that if a bipartite graph is antimagic, then
it has an antimagic orientation obtained by
directing all edges from one partite set to the other.
Thus by the result of Cranston~\cite{MR2478227},
regular bipartite graphs have an antimagic orientation. A bipartite graph is {\it biregular\/} if
vertices in each of  the same partite set have the same degree.  In this paper, by supporting
Conjecture~\ref{antimagic-orientation}, we obtain the
result below.

\begin{theorem}\label{th1}
Every biregular bipartite graph admits an antimagic orientation.
\end{theorem}

\section{Notation and Lemmas}


Let $G$ be a graph. If $G$ is bipartite with partite sets $X$ and $Y$,
we denote $G$ by $G[X,Y]$.
Given an orientation of $G$ and a labeling on $E(G)$,
for a vertex $v\in V(G)$ and a subgraph $H$ of $G$, we use $\omega_H(v)$ to denote
the {\it oriented sum at $v$\/} in $H$, which  is the sum of labels of all arcs entering $v$ minus
the sum of labels of all arcs leaving it in the graph $H$.
If $v$ is of degree 2 in $G$, we say the
labels at edges incident to $v$ the {\it label at $v$\/} and
write it as a pair in  $\{(a,b),(-a,b),(a,-b),(-a,-b)\}$, where $a,b$ are the labels on
the two edges incident to $v$, and $-a$ is used if
the edge with label $a$ is leaving $v$ and $a$ is used
otherwise; similar situation for the value $-b$ or $b$.

A \emph{trail} is an alternating sequence of vertices and edges $v_0,e_1,v_1,\ldots,e_t,v_t$ such that
$v_{i-1}$ and $v_i$ are the endvertices of $e_i$, for each $i$ with $1\leq i\leq t$, and the edges are all distinct (but there might be repetitions among the vertices).  A trail is \emph{open} if $v_0\neq v_t$. The \emph{length} of a trail is the number of edges in it.
Occasionally,  a trail $T$ is also treated as a graph whose
vertex set is the set of distinct vertices in $T$ and edge set is the set of edges in $T$.
We use the terminology ``trail'' without distinguishing if it is a sequence or a graph, but the meaning
will be clear from the context.
For two integers $a,b$ with $a<b$,  let $[a,b]:=\{a,a+1,\cdots, b\}$.

We need the result below which guarantees a matching
in a bipartite graph. A simple proof of this
result can be found in~\cite{local-CE}.

%

\begin{lemma} [\cite{local-CE}]\label{matching}
Let $H$ be a bipartite graph with partite sets
$X$ and $Y$.
If there is no isolated vertex in $X$
and $d_H(x)\ge d_H(y)$ holds
for every edge $xy$ with $x\in X$ and $y\in Y$,
then $H$ has a matching which saturates $X$.
\end{lemma}

For even regular graphs, Petersen proved that
a 2-factor always exists.

\begin{lemma}[\cite{MR1554815}]\label{petersen}
Every regular (multi)graph with positive even degree has a $2$-factor.
\end{lemma}

Also we need the following  result on
 decomposing edges in a graph into trails.

\begin{lemma}[\cite{MR3414180}]\label{trail}
Given a connected graph $G$, and let $T=\{v\in V: d_G(v)$ is odd$\}$.
If $T\neq \emptyset$, then $E(G)$ can be partitioned into $\frac{|T|}{2}$ open trails.
\end{lemma}

\begin{lemma}\label{cycle1}
Every simple $2$-regular graph $G$ admits a vertex distinguishable labeling with labels in $[a,b]$, where $a, b$ are two positive integers with $b-a=|E(G)|-1$. Moreover, the vertex sums belong to $[2a+1,2b-1]$.
\end{lemma}


\pf
Note that $G$ is antimagic by Corollary 3 in~\cite{MR2478227}.
Assume that $\phi: E(G)\rightarrow [1, |E(G)|]$ is an antimagic labeling of $G$.
Define another labeling $\varphi: E(G)\rightarrow [a,b]$ based on $\phi$ as follows.
\begin{align}
\varphi(e)=\phi(e)+a-1, \quad \forall \, e\in E(G).\nonumber
\end{align}
Since $G$ is regular and $\phi$ is antimagic, it is clear that $\varphi$
is a vertex distinguishable labeling of $G$.
Furthermore, the sums fall into the interval $[2a+1,2b-1]$.
\qed

\begin{lemma}\label{paths}
Let $T[X,Y]$ be an open trail
with all vertices in $Y$ having degree 2
except precisely two having degree 1.
Suppose $T$ has $2m$ edges.
Let $y_1$ and $y_{m+1}$ be the two degree 1 vertices in $Y$
such that $T$ starts at $y_1$ and ends at $y_{m+1}$.
Let $a, b$ be two integers with $a=b-2m+1$.
Then there exists a bijection from $ E(T)$ to $[a,a+m-1]\cup [b-m+1,b]=[a,b]$ such that
each of the following holds.
\begin{enumerate}[(i)]
  \item $\omega_T(x)=\frac{d_T(x)(a+b)}{2}$ for any $x\in X$; and
  $\omega_T(y)\ne \omega_T(z)$ for any distinct $y,z\in Y-\{y_1, y_{m+1}\}$.

  \item If  $m\equiv0 \,(\md 2)$,  then $\omega_T(y_1)=b$, $\omega_T(y_{m+1})=b-m+1$,
and  $\omega_T(y)$ is an odd number in $[2a+1,2a+2m-3]\cup [2b-2m+5,2b-3]$ for any $y\in Y-\{y_1, y_{m+1}\}$.

  \item If  $m\equiv 1 \,(\md 2)$,  then $\omega_T(y_1)=a$, $\omega_T(y_{m+1})=b-m+1$,
and
$\omega_T(y)$ is an odd number in $[2a+3,2a+2m-3]\cup [2b-2m+5,2b-1]$
for any  $y\in Y-\{y_1, y_{m+1}\}$.
\end{enumerate}
\end{lemma}

\pf
Since $|E(T)|=2m$, and except precisely two degree
1 vertices,
all other vertices in $Y$ have degree 2,
we conclude that $|Y|=m+1$.  Let $Y=\{y_1,y_2,\cdots, y_{m+1}\}$.
Then there are precisely $m$ edges of $T$ incident
to vertices in $Y$ with even indices,
and $m$ edges of $T$ incident
to vertices in $Y$ with odd indices. We  treat $T$
as an alternating sequence of vertices and edges starting at
$y_1$ and ending at $y_{m+1}$.

If $m\equiv0 \,(\md 2)$,
following the order of the appearances of
edges in $T$,
assign  edges incident to
vertices in $Y$ of even indices  with labels
$$
a, a+1, \cdots, a+m-1,
$$
and assign
edges incident to
vertices in $Y$ of odd indices  with labels
$$
b, b-1, \cdots, b-m+1.
$$
That is, the label at $y_i$ is $(a+i-2, a+i-1)$ if $i$ is even; and $(b-i+2, b-i+1)$ if $i$ is odd and not equal to $1$ or $m+1$.

If $m\equiv1 \,(\md 2)$,
following the order of the appearances of
edges in $T$,
assign  edges incident to
vertices in $Y$ of odd indices  with labels
$$
a, a+1, \cdots, a+m-1,
$$
and assign
edges incident to
vertices in $Y$ of even indices  with labels
$$
b, b-1, \cdots, b-m+1.
$$
That is, the label at $y_i$ is $(a+i-2, a+i-1)$ if $i$ is  odd and not equal to $1$ or $m+1$; and $(b-i+2, b-i+1)$ if $i$ is even.

If $m\equiv0 \,(\md 2)$, for each $y_i\in Y$ with $1\le i\le m+1$, by the assignment of labels, we have that
$$ \omega_T(y_i)=
\left\{
  \begin{array}{ll}
    b , & \hbox{if $i=1$;} \\
    b-m+1 , & \hbox{if $i=m+1$;} \\
    2a+2i-3, & \hbox{if $i$ is even and $2\le i\le m$;} \\
    2b-2i+3, & \hbox{if $i$ is odd and $3\le i\le m-1$.}
  \end{array}
\right.
$$
If $m\equiv1 \,(\md 2)$, for each $y_i\in Y$ with $1\le i\le m+1$, by the assignment of labels, we have that
$$ \omega_T(y_i)=
\left\{
  \begin{array}{ll}
    a , & \hbox{if $i=1$;} \\
    b-m+1 , & \hbox{if $i=m+1$;} \\
    2a+2i-3, & \hbox{if $i$ is odd and $3\le i\le m$;} \\
    2b-2i+3, & \hbox{if $i$ is even and $2\le i\le m-1$.}
  \end{array}
\right.
$$
The sum on each vertex $y_i$ with $y_i\in Y-\{y_1, y_{m+1}\}$ is expressed as either $2a+2i-3$ or $2b-2i+3$, which is an odd
number. Furthermore,
the sums on $y_2, y_4,\cdots, y_{m}$,  starting at $2a+1$, strictly increase to $2a+2m-3$ if $m\equiv0 \,(\md 2)$,
and  the sums on $y_3, y_5,\cdots, y_{m-1}$,  starting at $2a+3$, strictly increase to $2a+2m-3$ if $m\equiv1 \,(\md 2)$.
The sums on $y_3, y_5,\cdots, y_{m-1}$,  starting at $2b-3$,  strictly decrease  to $2b-2m+5$  if $m\equiv0 \,(\md 2)$,
and the sums on $y_2, y_4,\cdots, y_{m}$,  starting at $2b-1$,  strictly decrease  to $2b-2m+5$  if $m\equiv1 \,(\md 2)$.
So  these sums are all distinct.
Since $a=b-2m+1$,  it holds that $2a+2m-3<2b-2m+5$. Thus  all $\omega_T(y)$
are distinct for $y\in Y$  with $y\in Y-\{y_1, y_{m+1}\}$.

Let $x$ be a vertex in $X$. Suppose that one appearance of $x$ is adjacent
to $y_i$ and $y_{i+1}$ in the sequence $T$.
If $m\equiv0 \,(\md 2)$,
for  even $i$ with $2\le i\le m$, the labels on the two edges $xy_i$ and $xy_{i+1}$
contribute a value of $(a+i-1)+(b-(i+1)+2)=a+b$ to $\omega_T(x)$;
for odd $i$ with $1\le i\le m-1$, the labels on the two edges $xy_i$ and $xy_{i+1}$
contribute a value of $(b-i+1)+(a+(i+1)-2)=a+b$ to $\omega_T(x)$.
Since $x$ appears $d_T(x)/2$ times in $T$, $\omega_T(x)=\frac{d_T(x)(a+b)}{2}$.
If $m\equiv1 \,(\md 2)$,
for even $i$ with $2\le i\le m$, the labels on the two edges $xy_i$ and $xy_{i+1}$
contribute a value of $(b-i+1)+(a+(i+1)-2)=a+b$ to $\omega_T(x)$;
for odd $i$ with $1\le i\le m-1$, the labels on the two edges $xy_i$ and $xy_{i+1}$
contribute a value of $(a+i-1)+(b-(i+1)+2)=a+b$ to $\omega_T(x)$.
Since $x$ appears $d_T(x)/2$ times in $T$, $\omega_T(x)=\frac{d_T(x)(a+b)}{2}$.
\qed

%

\begin{lemma}\label{cycle2}
Let $C[X,Y]$ be a cycle of length $2m$ with
$m\equiv 0\,(\md 2)$, and
let $a, b$ be two integers with $a=b-2m+1$.
Then there exists a bijection from $ E(C)$ to $[a,a+m-1]\cup [b-m+1,b]=[a,b]$ such that
each of the following holds.
\begin{enumerate}[(i)]
  \item $\omega_C(x)=a+b$ for any $x\in X$.
  \item $\omega_C(y)\ne \omega_C(z)$ for any distinct $y,z\in Y$.
  \item $\omega_C(y)\in [2a+1,2a+2m-3]\cup [2b-2m+5,2b-2]$ for all $y\in Y$, and the sums in $[2a+1,2a+2m-3]$ are odd.
\end{enumerate}
\end{lemma}

\pf
Denote by  $C=x_1y_1x_2y_2\cdots x_my_mx_1$ with $x_i\in X$ and $y_i\in Y$.
Following the order of the appearances of
edges in $C$,
assign  edges incident to
vertices in $\{y_1, y_3,\cdots, y_{m-1}\}$  with labels
$$
a+1, a, a+2, a+3,\cdots, a+m-2, a+m-1.
$$
Note that the labels are increasing consecutive
integers after exchanging the positions of the
first two;
assign  edges incident to
vertices in $\{y_{2}, y_{4},\cdots, y_{m}\}$  with labels
$$
b,  b-2, b-3, \cdots,  b-m+2, b-m+1, b-1.
$$
Note that the labels are decreasing  consecutive
integers after inserting the last number between
the first two labels.

For each $y_i\in Y$ with $1\le i\le m$,
following the appearances of the
edges in the sequence of $C$,   we denote
the labels on the edges incident
to $y_i$ by an ordered pair.
Particularly, by the assignment of the labels, we have that
$$ \mbox{label at $y_i$}=
\left\{
  \begin{array}{ll}
    (a+1, a) , & \hbox{if $i=1$;} \\
    (b,b-2), & \hbox{if $i=2$;} \\
   (b-m+1, b-1) , & \hbox{if $i=m$;} \\
    (a+i-1,a+i), & \hbox{if $i$ is odd with $3\le i\le m-1$;}\\
    (b-i+1,b-i), & \hbox{if $i$ is even with  $4\le i\le m-2$.}
  \end{array}
\right.
$$
The sums on $y_1, y_3,\cdots, y_{m-1}$ starting at $2a+1$
strictly increase to $2a+2m-3$ and all of them are odd;
and the sums on $y_4, y_6,\cdots, y_{m-2}$ starting
at $2b-7$ strictly decrease to $2b-2m+5$ and all of them are odd; the sums on $y_2, y_m$
are even numbers $2b-2$ and $2b-m$, respectively.
Since $a= b-2m+1$, $2a+2m-3<2b-2m+5$. Hence  all $\omega_C(y)$
are distinct for $y\in Y$.
This shows both (ii) and (iii).

Let $x_i$ be a vertex in $X$ for $1\leq i\leq m$.
If $i=1$, the labels on the two edges $x_1y_1$ and $x_1y_{m}$
are $a+1$ and $b-1$, respectively;
if $i=2$, the labels on the two edges $x_2y_1$ and $x_2y_{2}$
are $a$ and $b$, respectively.  Thus $\omega_C(x_i)=a+b$
if $i=1,2$.
Suppose $i\ge 3$.
If $i$ is even, then the labels on the two edges $x_iy_{i-1}$ and $x_iy_i$
are $a+(i-1)$ and $b-i+1$, respectively;
if $i$ is odd, the labels on the two edges $xy_{i-1}$ and $xy_{i}$
are $b-(i-1)$ and $a+i-1$, respectively.
Thus $\omega_C(x_i)=a+b$.
This proves (i).
\qed


\begin{lemma}\label{pathcycle}
Let $G$ be a bipartite graph, and $H[X,Y]$ a subgraph of $G$.
Suppose
that $E(H)$ can be decomposed into edge-disjoint
$p+q$ open trails $T_1,  \cdots, T_p$, $T_{p+1}, \cdots, T_{p+q}$,
and $\ell$ cycles $C_{p+q+1}, \cdots, C_{p+q+\ell}$.
Suppose further that these $p+q+\ell$ subgraphs have no
common vertex in $Y$, and for each of the trail, its vertices in
$Y$ are all distinct and its  endvertices  are contained in $Y$.
Let $2m:=|E(H)|$, and $c, d$ be two integers with $c=d-2m+1$.
If the length of $T_1,\cdots, T_p$ are congruent to
2 modulo 4, and the length of each of the remaining  trails and cycles is
congruent to 0 modulo 4, then
there exists a bijection  from $ E(H)$ to $[c,c+m-1]\cup [d-m+1,d]=[c,d]$  such that
each of the following holds.
\begin{enumerate}[(i)]
\item  $\omega_H(x)=\frac{d_H(x)(c+d)}{2}$ for any $x\in X$.
 \item For each $i$ with $1\le i\le p+q$, let $T_i$ start at $y_{1i}$ and end at $y_{(m_i+1)i}$,
 where $y_{1i}, y_{(m_i+1)i}\in Y$. Suppose that
 $$
 \begin{cases}
   \omega_G(y_{1i})=\omega_H(y_{1i})+(c-2p+i-1),  & \hbox{if $1\le i\le p$}, \\
   \omega_G(y_{(m_i+1)i})=\omega_H(y_{(m_i+1)i})+(c-i),   & \hbox{if $1\le i\le p$},\\
      \omega_G(y_{1i})=\omega_H(y_{1i})+(d+2q-2(i-p-1)) ,& \hbox{if $p+1\le i\le p+q$},\\
      \omega_G(y_{(m_i+1)i})=\omega_H(y_{(m_i+1)i})+(d+2q-2(i-p-1)-1), & \hbox{if $p+1\le i\le p+q$},\\
      \omega_G(y)=\omega_H(y), & \text{if $y\ne y_{1i},y_{(m_i+1)i} $}.
\end{cases}
$$
Then $\omega_G(y)\ne \omega_G(z)$ for any distinct $y,z\in Y$,
and $\omega_G(y)\in [2c-2p,\max\{2d+2q-\sum\limits_{i=1}^pm_i, 2d\}]$ for all $y\in Y$.
\end{enumerate}
\end{lemma}

\pf
Since  $T_1,  \cdots, T_{p+q}$
and  $C_{p+q+1},  \cdots, C_{p+q+\ell}$
are edge-disjoint and pairwise have no common vertex in $Y$,
and for each of the trail, its vertices in
$Y$ are all distinct and its  endvertices  are contained in $Y$,
we conclude that in the graph $H$,
 all vertices in $X$ have even degree,
 all the endvertices of $T_i$
are precisely the degree 1 vertices in $Y$,
and  all other vertices in $Y$ have degree 2.
For each $i$ with $1\le i\le p+q$ and each $j$ with $p+q+1\le j\le p+q+\ell$,
let
$$
|E(T_i)|=2m_i, \quad  |E(C_j)|=2m_j.
$$
Since  $V(T_i)\cap Y$
contains two degree 1 vertices  and $|V(T_i)\cap Y|-2$
degree 2 vertices in $Y$, we conclude that $|V(T_i)\cap Y|=m_i+1$.

Apply Lemma~\ref{paths} on each $T_i$, $1\le i\le p+q$,
with
$$
a:=a_i:=c+\sum\limits_{j=1}^{i-1}m_j, \quad b:=b_i:=d-\sum\limits_{j=1}^{i-1}m_j;
$$
and apply Lemma~\ref{cycle2} on each $C_i$, $p+q+1\le i\le p+q+\ell$,
with
$$
a:=a_i:=c+\sum\limits_{j=1}^{i-1}m_j, \quad b:=b_i:=d-\sum\limits_{j=1}^{i-1}m_j.
$$
Note that $a_i+b_i=c+d$, $1\le i\le p+q+\ell$.
By Lemma~\ref{paths} and Lemma~\ref{cycle2}, we get that
$$
\omega_H(x)=\frac{d_H(x)(c+d)}{2} \quad \mbox{for any $x\in X$.}
$$
By Lemma~\ref{paths}, the sums at $y_{1i},y_{m_ii}$, respectively, are
\begin{numcases}{}
\omega_H(y_{1i})=c+\sum\limits_{j=1}^{i-1}m_j, \quad \omega_H(y_{(m_i+1)i})=d-\sum\limits_{j=1}^{i}m_j+1, & \mbox{ if $1\le i\le p$};\label{sumy1}\\
\omega_H(y_{1i})=d-\sum\limits_{j=1}^{i-1}m_j,  \quad \omega_H(y_{(m_i+1)i})=d-\sum\limits_{j=1}^{i}m_j+1,& \mbox{ if $p+1\le i\le p+q$};\label{sumy2}
\end{numcases}
and for each $i$ with $1\le i\le p+q$,  the  sums at
vertices in $V(T_i)\cap Y-\{y_{1i},  y_{(m_i+1)i}\}$
fall into the intervals
$$
\begin{cases}
    [2a_i+3, 2a_i+2m_i-3]\cup [2b_i-2m_i+5, 2b_i-1],  & \text{if $1\le i\le p$ }, \\
      [2a_i+1, 2a_i+2m_i-3]\cup [2b_i-2m_i+5, 2b_i-3],& \text{if $p+1\le i\le p+q$},
\end{cases}
$$
and all these sums are distinct and odd.

By Lemma~\ref{cycle2},
for each $i$ with $p+q+1\le i\le p+q+\ell$,  the  sums at
vertices in $V(C_i)\cap Y$ are all distinct and
fall into the intervals
$$
[2a_i+1, 2a_i+2m_i-3]\cup [2b_i-2m_i+5, 2b_i-2],
$$
and all the sums in $[2a_i+1, 2a_i+2m_i-3]$ are odd.
%

Since for each $i,j$ with $1\le i< j\le p+q+\ell$,
$$
\left\{
  \begin{array}{ll}
   2a_j+1>2a_i+2m_i-3;\\
   2b_i-2m_i+5>2b_j-1,
    \end{array}
\right.
$$
we see that
$2a_{p+q+\ell}+2m_{p+q+\ell}-3$ is the largest value in the set
$$
\left(\bigcup_{i=1}^{p}[2a_i+3, 2a_i+2m_i-3]\right)\bigcup \left(\bigcup_{i=p+1}^{p+q+\ell}[2a_i+1, 2a_i+2m_i-3]\right),
$$
and
$2b_{p+q+\ell}-2m_{p+q+\ell}+5$ is the smallest value in the set
\begin{eqnarray*}
\left(\bigcup_{i=1}^{p}[2b_i-2m_i+5, 2b_i-1]\right)&\bigcup& \left(\bigcup_{i=p+1}^{p+q}[2b_i-2m_i+5, 2b_i-3]\right)\\
& \bigcup& \left(\bigcup_{i=p+q+1}^{p+q+\ell}[2b_i-2m_i+5, 2b_i-2]\right).
\end{eqnarray*}
Furthermore,
\begin{eqnarray*}
\begin{tabular}{lll}
  &$2b_{p+q+\ell}-2m_{p+q+\ell}+5-(2a_{p+q+\ell}+2m_{p+q+\ell}-3)$ &\\
=& $2d-\sum\limits_{j=1}^{p+q+\ell}2m_i+5-(2c+\sum\limits_{j=1}^{p+q+\ell}2m_i-3)$&\\
  $=$& $2(c+2m-1)-2m+5-(2c+2m-3)$\quad \mbox{($d= c+2m-1$ and $\sum\limits_{j=1}^{p+q+\ell}2m_i=2m$)}& \\
   $=$& $6$.  &
\end{tabular}
\end{eqnarray*}
Hence, all the vertex  sums at $Y-\{y_{1i},  y_{(m_i+1)i}\,|\,  1\le i\le p+q\}$
are pairwise distinct.
Let
$$
k_1:=\sum\limits_{i=1}^p m_i, \quad \mbox{and}\quad  k_2:=\sum\limits_{i=p+1}^{p+q} m_i.
$$
By  Equalities (\ref{sumy1}) and (\ref{sumy2}),  and the assumptions on the parity of each $m_i$, if $1\le i\le p$,
 \begin{numcases}{}
   \omega_G(y_{1i})=\omega_H(y_{1i})+(c-2p+i-1),   \nonumber  \\
   \hphantom{\omega_G(y_{1i})}=c+\sum\limits_{j=1}^{i-1}m_j+c-2p+i-1=2c-2p+\sum\limits_{j=1}^{i-1}(m_j+1) \nonumber  \\
  \hphantom{\omega_G(y_{1i})} \equiv 0\,(\md 2) &\text{since each $m_j$ is odd }\nonumber \\
   \hphantom{\omega_G(y_{1i})}\in [2c-2p, 2c-p-1+k_1-m_p], \label{y1ip}\\
   \omega_G(y_{(m_i+1)i})=\omega_H(y_{(m_i+1)i})+(c-i)  \nonumber\\
   \hphantom{\omega_G(y_{(m_i+1)i})}=d-\sum\limits_{j=1}^{i}m_j+1+c-i\nonumber \\
   \hphantom{\omega_G(y_{(m_i+1)i})}=c+2m-1-\sum\limits_{j=1}^{i}m_j+1+c-i\nonumber\\
     \hphantom{\omega_G(y_{(m_i+1)i})}=2c+2m-\sum\limits_{j=1}^{i}(m_j-1)\equiv 0\,(\md 2)&\text{since each $m_j$ is odd }\nonumber\\
    \hphantom{\omega_G(y_{(m_i+1)i})}\in [2c+2m-k_1-p, 2c+2m-m_1-1], \label{ym1p}
    \end{numcases}
and if $p+1\le i\le p+q$,
\begin{numcases}{}
      \omega_G(y_{1i})=\omega_H(y_{1i})+(d+2q-2(i-p-1)) \label{y1iq}\\
     \hphantom{\omega_G(y_{1i})}= d-\sum\limits_{j=1}^{i-1}m_j+d+2q-2(i-p-1)\equiv 0\,(\md 2),&\text{since each $m_j$ is even }\nonumber\\
      \omega_G(y_{(m_i+1)i})=\omega_H(y_{(m_i+1)i})+(d+2q-2(i-p-1)-1) \label{ym1q}\\
      \hphantom{\omega_G(y_{(m_i+1)i})}= d-\sum\limits_{j=1}^{i}m_j+1+d+2q-2(i-p-1)-1\nonumber\\
      \hphantom{\omega_G(y_{(m_i+1)i})}\equiv 0\,(\md 2),&\text{since each $m_j$ is even }\nonumber
\end{numcases}

By the above analysis, for each $i$ with  $1\le i\le p+q$, both
$\omega_G(y_{1i})$   and $\omega_G(y_{(m_i+1)i})$
are even.
As all the sums at vertices in $\bigcup_{i=1}^{p+q} V(T_i)\cap Y-\{y_{1i},  y_{(m_i+1)i}\,|\,  1\le i\le p+q\}$
are  odd  by Lemma~\ref{paths}, and by Lemma~\ref{cycle2} all the sums  on vertices in
$\bigcup_{i=p+q+1}^{p+q+\ell}V(C_i)\cap Y$ which fall into the set  $\bigcup_{i=p+q+1}^{p+q+\ell}[2a_i+1, 2a_i+2m_i-3]$ are odd, all of them are distinct from these $2(p+q)$ $\omega_G$ sums on vertices in $\{y_{1i},  y_{(m_i+1)i}\,|\,  1\le i\le p+q\}$.
Hence, to show that all the vertex  sums at vertices in
$Y$ are distinct, we are left  to check that all these $2(p+q)$ $\omega_G$ sums
are distinct with the sums
on vertices in
$\bigcup_{i=p+q+1}^{p+q+\ell}V(C_i)\cap Y$ which fall into the set  $\bigcup_{i=p+q+1}^{p+q+\ell}[2b_i-2m_i+5, 2b_i-2]$,
and all these $2(p+q)$ $\omega_G$ sums at vertices in $\{y_{1i},  y_{(m_i+1)i}\,|\,  1\le i\le p+q\}$
are pairwise distinct.

If $1\le i<j\le p$, by~(\ref{y1ip}) and (\ref{ym1p}),
$\omega_G(y_{1i})<\omega_G(y_{1j})$ and $\omega_G(y_{(m_i+1)i})>\omega_G(y_{(m_i+1)j})$.
Thus, all the sums at vertices either  in $\{ y_{1i}\,|\, 1\le i\le p\}$
or in $\{ y_{(m_i+1)i}\,|\, 1\le i\le p\}$ are all distinct, and
$$
\omega_G(y_{1p})=\max\{\omega_G(y_{1i})\,|\, 1\le i\le p\}, \quad
\omega_G(y_{(m_{p}+1)p})=\min\{\omega_G(y_{(m_i+1)i})\,|\, 1\le i\le p\}.
$$
Furthermore, by~(\ref{y1ip}) and (\ref{ym1p}),
\begin{eqnarray*}
\begin{tabular}{lll}
  &$\omega_G(y_{(m_{p}+1)p})-\omega_G(y_{1p})$&\\
=& $2c+2m-k_1-p-(2c-p-1+k_1-m_p)$&\\
  =& $2m-2k_1+m_p+1\ge 2$\quad \mbox{($2m\ge 2k_1, m_p\ge 1$).}&
  \end{tabular}
\end{eqnarray*}
Thus, the $\omega_G$ sums on vertices in $\{ y_{1i}, y_{(m_i+1)i}\,|\, 1\le i\le p\}$ are all distinct.

 By~(\ref{y1iq}) and (\ref{ym1q}),
$\omega_G(y_{1i})>\omega_G(y_{(m_i+1)i})$  for all $i$ with $p+1\le i\le p+q$,
and $\omega_G(y_{(m_i+1)i})>\omega_G(y_{1j})$ if $p+1\le i<j\le p+q$.
Thus, all the $\omega_G$ sums on vertices in $\{ y_{1i}, y_{(m_i+1)i}\,|\, p+1\le i\le p+q\}$ are all distinct,
and
$$
\omega_G(y_{(m_{p+q}+1)(p+q)})=\min\{\omega_G(y_{1i}),\omega_G(y_{(m_i+1)i})\,|\, p+1\le i\le p+q\}.
$$
Furthermore,
$$
\omega_G(y_{(m_{1}+1)1})=\max\{\omega_G(y_{1i}),\omega_G(y_{(m_i+1)i})\,|\, 1\le i\le p\},
$$
and by~(\ref{ym1q}) and (\ref{ym1p}),
\begin{eqnarray*}
\begin{tabular}{lll}
  &$\omega_G(y_{(m_{p+q}+1)(p+q)})-\omega_G(y_{(m_1+1)1})$&\\
=& $2d-k_1-k_2+2-(2c+2m-m_1-1)$&\\
  =& $2m-k_1-k_2+m_1+1\ge 2$\quad \mbox{($d=c+2m-1, 2m\ge k_1+k_2, m_1\ge 1$).}&
  \end{tabular}
\end{eqnarray*}
Thus, all the $\omega_G$ sums on vertices in $\{ y_{1i}, y_{(m_i+1)i}\,|\, 1\le i\le p+q\}$ are all distinct.
We may assume that $\ell\ge 1$. Otherwise, we are done.

By the definition of the parameters $b_i$ and easy calculations,
$$
\bigcup_{i=p+q+1}^{p+q+\ell}[2b_i-2m_i+5, 2b_i-2]\subseteq [2b_{p+q+\ell}-2m_{p+q+\ell}+5, 2b_{p+q+1}-2].
$$
By~(\ref{y1ip}), (\ref{ym1p}),(\ref{y1iq}), and (\ref{ym1q}), the
$\omega_G$ sums at vertices in $\{ y_{1i}, y_{(m_i+1)i}\,|\, 1\le i\le p+q\}$
fall into the intervals
$$
[2c-2p, 2c+2m-m_1+1]\cup [2d-k_1-k_2+2, 2d+2q-k_1].
$$
Since
\begin{eqnarray*}
\begin{tabular}{lll}
  &$(2b_{p+q+\ell}-2m_{p+q+\ell}+5)-(2c+2m-m_1-1)$ &\\
=& $2d-2k_1-2k_2+5-(2(d-2m+1)+2m-m_1-1)$&\\
  =& $2m-2k_1-2k_2+m_1+4$& \\
   $\ge $& $5$  \quad \mbox{($2m\ge 2k_1+2k_2$ and $m_1\ge 1$)},&
\end{tabular}
\end{eqnarray*}
and
\begin{eqnarray*}
\begin{tabular}{lll}
  &$(2d-k_1-k_2+2)-(2b_{p+q+1}-2)$ &\\
=& $2d-2k_1-k_2+2-(2d-2k_1-2k_2-2)$&\\
  =& $k_1+k_2+4 \ge 4$,&
  \end{tabular}
\end{eqnarray*}
%
we then conclude that these $2(p+q)$ $\omega_G$ sums at vertices in $\{y_{1i},  y_{(m_i+1)i}\,|\,  1\le i\le p+q\}$
are all distinct with  the $\omega_G$ sums at vertices in $\bigcup_{i=p+q+1}^{p+q+\ell}V(C_i)\cap Y$ which fall into
the set $\bigcup_{i=p+q+1}^{p+q+\ell}[2b_i-2m_i+5, 2b_i-2]$.

By all the  arguments above, we have shown that
$\omega_G(y)\ne \omega_G(z)$ for any distinct $y,z\in Y$.
Since the sums on vertices in $Y-\{y_{1i},  y_{(m_i+1)i}\,|\,  1\le i\le p+q\}$
fall into the interval $[2c, 2d]$ and the sums on vertices in
$\{y_{1i},  y_{(m_i+1)i}\,|\,  1\le i\le p+q\}$ fall into the interval  $[2c-2p,2d+2q-k_1]$,
where the value $2d+2q-k_1$ is attained at $\omega_G(y_{1(p+1)})$,
we have that  $\omega_G(y)\in [2c-2p,\max\{2d+2q-k_1, 2d\}]$ for all $y\in Y$.
\qed

\section{Proof of Theorem \ref{th1}}
Let $G=[X,Y]$ be a biregular bipartite graph.
Assume that $|X|=m$, $|Y|=n$, $d_G(x)=s\geq d_G(y)=t$, where $x\in X, y\in Y$.
Consequently $m\leq n$ and $|E(G)|=ms=nt$.
Given an orientation of $G$, we will denote
the orientation  by $\overrightarrow{G}$.

If $t=1$,
then  $G$ is the union of vertex-disjoint stars
with centers in $X$.
Denote
$$
X=\{x_1,x_2,\cdots,x_m\} \quad \mbox{and} \quad
Y=\{y_1,y_2,\cdots, y_n\}.
$$
For each $x_i$, $1\le i\le m$,  we  assign arbitrarily
edges incident to $x_i$ with labels
$$
s(i-1)+1, s(i-1)+2,\cdots, s(i-1)+s.
$$

Orient edges of $G$ from $X$ to $Y$.
Thus, the oriented vertex sums for vertices in $X$
are negative,  and the oriented
vertex sums for vertices in $Y$
are positive. Hence,
no two vertices $x$ and $y$ conflict if $x\in X$ and
$y\in Y$.
Also, it is routine to check that
no two vertices in $X$ conflicting and
no two vertices in $Y$ conflicting.
Hence the labeling of $\overrightarrow{G}$ is
antimagic.
Thus we assume $t\ge 2$. We
distinguish three  cases for finishing
the proof.

%
%
%
\bigskip
{\noindent \bf  Case 1: $t\ge 3$}

Orient edges of $G$ from $X$ to $Y$, and denote the orientation
by $\overrightarrow{G}$. By the orientation of $G$, the sums of vertices in $X$ are negative while
the sums  at vertices in $Y$ are positive. Hence in the following,
we just need to find a labeling of $\overrightarrow{G}$ using labels in $[1,sm]$,
 which guarantees that the sums at vertices in $X$ are all distinct and the sums at vertices in $Y$ are all distinct.
By Lemma \ref{matching}, $G$ has a matching $M$ saturating $X$.
Assume, w.l.o.g,  that $M=\{x_1y_1, x_2y_2,\ldots,x_my_m\}$.
Let $H=G-M$. Note that $d_H(y_i)=t-1$ for $1\leq i\leq m$ and $d_H(y_i)=t$ for $m+1\leq i\leq n$.

\bigskip
{\noindent \bf  Subcase 1.1: $t\ge 3$ and $t$ is odd}

Reserve labels in $[1,m]$ for edges in $M$, and use labels in
$[m+1, tn=sm]$ for edges in $H$.
For each $y_i$ with $1\le i\le m$,
assign arbitrarily the edges incident to $y_i$
 with labels
$$
m+i, 3m-i+1, 3m+i, 5m-i+1, 5m+i, \cdots, (t-2)m-i+1, (t-2)m+i, tm;
$$
and for each $y_i$ with $m+1\le i\le n$,
assign arbitrarily the edges incident to $y_i$
with labels
$$
t(i-m-1)+tm+1, t(i-m-1)+tm+2, \cdots, t(i-m-1)+tm+t.
$$

Assume, w.l.o.g.,  that under the above assignment
of labels, $\omega_H(x_1)\le \omega_H(x_2)\le \cdots \le \omega_H(x_m)$.
Now for each edge $x_iy_i\in M$, $1\le i\le m$,  assign the edge $x_iy_i$ with the label
$i$.

We verify now that the labeling of $\overrightarrow{G}$ given above is antimagic.
For each $x_i, x_j\in X$ with $i<j$, since
$\omega_H(x_i)\le \omega_H(x_j)$, it holds that $\omega_G(x_i)=\omega_H(x_i)+i<\omega_H(x_j)+j=\omega_G(x_j)$.

Next for each $y_i, y_j\in \{y_1,y_2,\cdots, y_m\}$ with $i<j$,
since
$\omega_H(y_i)=\omega_H(y_j)=\frac{t-1}{2}((t+1)m+1)$,
we have that  $\omega_G(y_i)=\omega_H(y_i)+i<\omega_H(y_j)+j=\omega_G(y_j)$.
By the assignment of labels on edges incident to $y_i$ with $m+1\le i\le n$,
the sums at $y_i$ are pairwise distinct.
The smallest vertex sum among these values is $\omega_G(y_{m+1})=t^2m+\sum\limits_{i=1}^{t}i$.
The largest vertex sum among values in $\{\omega_G(y_1),\cdots, \omega_G(y_m)\}$
is $\omega_G(y_m)=\frac{t-1}{2}((t+1)m+1)+m$.  It is easy to check that
$\omega_G(y_m)<\omega_G(y_{m+1})$. Hence, all the sums at vertices in $Y$
are distinct.


\bigskip
{\noindent \bf  Subcase 1.2: $t\ge 3$ and $t$ is even}

Reserve labels in $\{2,4,\cdots, 2m\}$ for edges in $M$,
and use the labels in $\{1,3,\cdots, 2m-1\}\cup \{2m+1,\cdots, tn=sm\}$
for edges in $H$.
For each $y_i$ with $1\le i\le m$,
assign arbitrarily the edges incident to $y_i$
 with labels
$$
2i-1, 3m-i+1, 4m-i+1, 4m+i, 6m-i+1, \cdots, (t-2)m+i, tm-i+1;
$$
and for  each $y_i$ with $m+1\le i\le n$,
assign arbitrarily the edges incident to $y_i$
with labels
$$
t(i-m-1)+tm+1, t(i-m-1)+tm+2, \cdots, t(i-m-1)+tm+t.
$$
Assume, w.l.o.g.,  that under the above assignment
of labels, $\omega_H(x_1)\le \omega_H(x_2)\le \cdots \le \omega_H(x_m)$.
Now for each edge $x_iy_i\in M$, $1\le i\le m$,  assign the edge $x_iy_i$ with the label
$2i$.

We verify now that the labeling of $\overrightarrow{G}$ given above is antimagic.
Obviously, for each $x_i, x_j\in X$ with $1\leq i<j\leq m$,
because
$\omega_H(x_i)\le \omega_H(x_j)$, it holds that $\omega_G(x_i)=\omega_H(x_i)+2i<\omega_H(x_j)+2j=\omega_G(x_j)$.

Next for each $y_i, y_j$ with $1\leq i<j\leq m$, since
$\omega_H(y_i)=\omega_H(y_j)=(t^2-9)m+t-3$, we have that $\omega_G(y_i)=\omega_H(y_i)+2i<\omega_H(y_j)+2j=\omega_G(y_j)$.
By the assignment of labels on edges incident to $y_i$ with $m+1\le i\le n$,
the sums at $y_i$ are pairwise distinct.
The smallest sum among these values is $\omega_G(y_{m+1})=t^2m+\sum\limits_{i=1}^{t}i$.
The largest sum among values in $\{\omega_G(y_1),\cdots, \omega_G(y_m)\}$
is $\omega_G(y_m)=(t^2-9)m+t-3+2m$.  It is easy to check that
$\omega_G(y_m)<\omega_G(y_{m+1})$. Hence, all the sums at vertices in $Y$
are distinct.

\bigskip

{\noindent \bf  Case 2: $t=2$ and $s$ is odd}

By Lemma \ref{matching}, there exists a matching $M$ saturating vertices in $X$.
In each component of $G-M$, the vertices contained in $X$ are all of even degree $s-1$,
and all vertices contained in $Y$ are of degree 2 or 1.
Thus, the number of vertices with  degree 1 in the component is even.
Since there are in total $m$ vertices of degree 1 in $Y$,
by Lemma \ref{trail}, we can decompose $E(G-M)$ into $m/2$ open trials
with endvertices in $Y$.
Denote the trails by $T_1, \cdots, T_{m/2}$.
Since the endvertices of each $T_i$ are in $Y$,
all the vertices in $V(T_i)\cap X$ have even degree.
Consequently, $T_i$ has even length.
For
 each  $i$ with $1\le i\le m/2$,
let
$$
2m_i:=|E(T_i)|  \quad \mbox{and}\quad  T_i=y_{1i}x_{1i}\cdots x_{m_ii}y_{(m_i+1)i},
$$
where $x_{1i},x_{2i},\cdots, x_{m_ii}\in X$
and $y_{1i},y_{2i},\cdots, y_{m_i+1i}\in Y$. Note that
$x_{1i},x_{2i},\cdots, x_{m_ii}$ maynot be distinct
vertices in $X$, but $y_{1i},y_{2i},\cdots, y_{m_i+1i}$
are distinct vertices in $Y$ because  all vertices in $Y$
have degree 2 in $G$.
Assume further, w.l.o.g.,  that there are $p$ trails $T_1,\cdots, T_p$
of length congruent to 2 modulo 4,
and $q$  trails   $T_{p+1},\cdots, T_{p+q}$ of length  congruent to 0 modulo 4.
Set
$$
c:=2p+1,\quad  d:=sm-2q, \quad \text{and} \quad H=\bigcup_{i=1}^{p+q} T_i.
$$
The endvertices of the $m/2$ open trails
are exactly the set of Y-endvertices of the
$m$ matching edges.
Thus,  for each $i$ with $1\le i\le m/2$,
and for each edge  $e\in M$,  $e$
is incident to either  $y_{1i}$ or $y_{(m_i+1)i}$.
We assign
labels in $[1,2p]\cup [sm-2q+1, sm]$ on $e$ as below.
If $1\le i\le p$,
\begin{numcases}{\mbox{label on $e=$}}
     i, & \hbox{if $e$ is incident to $y_{1i}$;} \label{m1}\\
   2p-i+1, & \hbox{if $e$ is incident to $y_{(m_i+1)i}$;} \label{m2}
  \end{numcases}
 if $p+1\le i\le p+q$,
\begin{numcases}{\mbox{label on $e=$}}
    sm-2(i-p-1), & \hbox{if $e$ is incident to $y_{1i}$;} \label{m3}\\
   sm-2(i-p-1)-1, & \hbox{if $e$ is incident to $y_{(m_i+1)i}$.} \label{m4}
  \end{numcases}
Thus,
$$
 \begin{cases}
   \omega_G(y_{1i})=\omega_H(y_{1i})+i,  & \text{if $1\le i\le p$}, \\
   \hphantom{\omega_G(y_{1i})}=\omega_H(y_{1i})+(c-2p+i-1),\\
   \omega_G(y_{(m_i+1)i})=\omega_H(y_{(m_i+1)i})+(2p-i+1),   & \text{if $1\le i\le p$},\\
   \hphantom{\omega_G(y_{(m_i+1)i})}=\omega_H(y_{(m_i+1)i})+(c-i),\\
      \omega_G(y_{1i})=\omega_H(y_{1i})+(sm-2(i-p-1)) ,& \text{if $p+1\le i\le p+q$},\\
      \hphantom{\omega_G(y_{1i})}=\omega_H(y_{1i})+(d+2q-2(i-p-1)),\\
      \omega_G(y_{(m_i+1)i})=\omega_H(y_{(m_i+1)i})+(sm-2(i-p-1)-1), & \text{if $p+1\le i\le p+q$},\\
      \hphantom{\omega_G(y_{(m_i+1)i})}=\omega_H(y_{(m_i+1)i})+(d+2q-2(i-p-1)-1),\\
      \omega_G(y)=\omega_H(y), & \text{if $y\in Y, y\ne y_{1i},y_{(m_i+1)i} $}.
\end{cases}
$$

Applying  Lemma~\ref{pathcycle} on $H$ with $c=2p+1$ and $d=sm-2q$
defined as above,
we get an assignment of  labels  on $E(H)$
such that
\begin{enumerate}[(i)]
 \item For any $x\in V(H)\cap X$, $\omega_H(x)=\frac{d_H(x)(c+d)}{2}=\frac{(s-1)(sm-2q+2p+1)}{2}$; and
 \item For any distinct $y, z\in V(H)\cap Y$, $\omega_{G}(y)\ne \omega_G(z)$.
\end{enumerate}

Orient all the edges
of $G$ from $X$ to $Y$,
and denote the  orientation  by
$\overrightarrow{G}$.

\bigskip
{\noindent \bf Claim 1:}  The labeling of $\overrightarrow{G}$
given above  is antimagic.

\pf We first show that all the set of labels used is the set $[1,sm]$.
The set of labels used on edges in $M$
is the set $[1,2p]\cup [sm-2q+1,sm]$,
and the
set of labels used on edges in $H$ is the set $[c,d]=[2p+1, sm-2q]$.
The union of these two sets is the set  $[1,sm]$.

We then show that all the oriented sums on
vertices  in $\overrightarrow{G}$ are pairwise distinct.
We first examine the sums on vertices in $Y$.
For any
$y\in Y$,  since $Y=V(H)\cap Y$,
we know that all the sums on vertices in $Y$ are pairwise distinct
by (ii) preceding  Claim 1.

Next we show that in $\overrightarrow{G}$,  the oriented sums on
vertices in $X$ are all distinct.
According to (i) preceding Claim 1 and the orientation of $G$, we see that  for any $x\in X$,
$$
\omega_H(x)=\frac{-(s-1)(sm-2q+2p+1)}{2}.
$$
Since all the labels on the edges in $M$ are
distinct, and for any $x\in X$,
$$
\omega_G(x)=\omega_H(x)-\text{the label on  $e\in M$ which is incident to $x$},
$$
we know that the oriented sums on
vertices in $X$ are all distinct.


Finally, we show that for any $x\in X$ and $y\in Y$
the oriented sums at $x$ and $y$ are distinct in $G$.
This is clear since all the oriented sums at vertices in $Y$
are positive while that at vertices in $X$ are negative.
\qed

\bigskip

{\noindent \bf  Case 3: $t=2$ and $s$ is even}

We may assume that $s\ge 4$.  Otherwise $G$ is
2-regular and $|E(G)|=2m$. By Lemma~\ref{cycle1}, $G$ has an antimatic labeling
by taking $a:=1$ and $b:=2m$, and the labeling is
also an antimagic labeling of $\overrightarrow{G}$ obtained by
orienting all edges from $X$ to $Y$.

\bigskip
{\noindent \bf Claim 2:} The graph $G[X,Y]$ contains a
subgraph $F$ such that
\begin{itemize}
  \item [$(1)$] $F$ is a set of vertex disjoint cycles; and
  \item [$(2)$]$V(F)\cap X=X$.
\end{itemize}

\pf Suppressing all degree 2 vertices in $Y$, we obtain
an $s$-regular (multi)graph $G'$.
Since $s$ is even, by applying
Lemma~\ref{petersen}, we find a 2-factor of $G'$.
Subdivide each edge
in the 2-factor of $G'$, we get
the desired graph $F$.
\qed

Now $G-E(F)$ is a graph with all vertices having even degree.
So $G-E(F)$ can be decomposed into edge-disjoint cycles.
Assume that there are in total $\ell$ edge-disjoint cycles in
$G-E(F)$ such that
each of them has length congruent to 0 modulo 4,
and  there are in total $h$ edge-disjoint cycles in
$G-E(F)$ such that
each of them has length congruent to 2 modulo 4.
For each $i$ with $1\le i\le h$, denote  by
$$
C_i=x_{1i}y_{1i}\cdots x_{m_i i}y_{m_i i}x_{1i}, \quad \mbox{where}~ x_{1i},x_{2i}, \cdots, x_{m_i i}\in X,\quad y_{1i},y_{2i}, \cdots, y_{m_i i}\in X.
$$
the $i$-th cycle of length congruent to 2 modulo 4.

We pre-label edges in $\{x_{1i}y_{1i}, x_{1i}y_{m_ii}, y_{1i}x_{2i}, y_{m_ii}x_{m_ii}, x_{2i}y_{2i}, x_{m_ii}y_{(m_i-1)i}\,|\,  1 \le i\le h\}$.
In doing so, we distinguish if $s=4$ or $s\ge 6$.

If $s=4$, for each $i$ with $1\le i\le h$,  use the labels in $[1,3h]\cup [2m-3h+1,2m]$ to label
each edge $e$ indicated below.
\begin{numcases}
{\mbox{label on $e =$}}
      \ i,  & \hbox{if $e=x_{1i}y_{1i}$;} \label{s43}\\
   \ 2m-(i-1), & \hbox{if $e=x_{1i}y_{m_ii}$;} \label{s44}\\
   \ h+2i-1,  & \hbox{if $e=y_{1i}x_{2i}$;} \label{s45} \\
   \ h+2i,   &\hbox{if $e=y_{m_ii}x_{m_ii}$;}\label{s46}\\
   \ 2m+1-(h+2i-1),  &\hbox{if $e=x_{2i}y_{2i}$;} \label{s41}\\
    \ 2m+1-(h+2i),  & \hbox{if $e=x_{m_ii}y_{(m_i-1)i}$.\label{s42}}
  \end{numcases}
If $s\ge 6$, for each $i$ with $1\le i\le h$, use the labels in $[1,4h]\cup [(s-2)m-2h+1,(s-2)m]$ to label each edge $e$ indicated below.
\begin{numcases}
{\mbox{label on $e =$}}
     i, & \hbox{if $e=x_{1i}y_{1i}$;} \label{s63}\\
   h+i, & \hbox{if $e=x_{1i}y_{m_ii}$;}\label{s64} \\
    2h+2i-1, & \hbox{if $e=y_{1i}x_{2i}$;}\label{s65} \\
    2h+2i, & \hbox{if $e=y_{m_ii}x_{m_ii}$;}\label{s66}\\
   (s-2)m+2h+1-(2h+2i-1), & \hbox{if $e=x_{2i}y_{2i}$;}\label{s61} \\
    (s-2)m+2h+1-(2h+2i), & \hbox{if $e=x_{m_ii}y_{(m_i-1)i}$.} \label{s62}
  \end{numcases}

Assume that there are $q$ paths with positive length
after deleting the
vertices $x_{1i}, x_{2i}, x_{m_ii}$, $y_{1i}, y_{m_ii}$ in each  $C_i$ for $1\leq i\leq h$.
Assume, w.l.o.g., that these paths are obtained from
$C_{h-q+1}, \cdots, C_h$. For each $1\le i\le q $, denote these paths
by
$$
P_i=C_{h-q+i}-\{x_{1(h-q+i)}, x_{2(h-q+i)}, x_{m_{(h-q+i)}(h-q+i)}, y_{1(h-q+i)}, y_{m_{(h-q+i)}(h-q+i)}\},
$$
and assume that $P_i$ starts at $y_{2(h-q+i)}$
and ends at $y_{(m_{h-q+i}-1)(h-q+i)}$.
Then each $P_i$ has length $2(m_{h-q+i}-3)\equiv 0(\md4)$.
Under this assumption, we know that $C_1, C_2, \cdots, C_{h-q} $
are 6-cycles.

Denote the  $\ell$ edge-disjoint cycles in
$G-E(F)$ such that
each of them has length congruent to 0 modulo 4 by $D_{q+1}, D_{q+2}, \cdots, D_{q+\ell}$.
Let
$$
H=\left(\bigcup_{i=1}^q P_i\right)\bigcup \left(\bigcup_{i=q+1}^{q+\ell} D_i\right).
$$

If $s=4$, then
let
$$
c:=3h+1 \quad \mbox{and} \quad d:=2m-3h.
$$
It is clear that $d-c=2m-6h-1=|E(H)|-1$.
For $h-q+1\le j\le h$, let $i=j-h+q$. Then
by Equations~(\ref{s41}) and (\ref{s42}),
the labels on the
other edges not in $H$ incident to
$y_{2j}$ and $y_{(m_j-1)j}$,  respectively,
are
\begin{eqnarray*}
2m+1-(h+2j-1)&=& (s-2)m-h-2(i+h-q)+2\\
                               &=&(s-2)m-3h+2q-2i+2=d+2q-2(i-1),\\
  2m+1-(h+2j)&=&(s-2)m-h-2(i+h-k)+1\\
                             &=&(s-2)m-3h+2q-2i+1=d+2q-2(i-1)-1.
\end{eqnarray*}
If $s\ge 6$, then
let
$$
c:=4h+1 \quad \mbox{and} \quad d:=(s-2)m-2h.
$$
Again $d-c=(s-2)m-6h-1=|E(H)|-1$.
For $h-q+1\le j\le p$, let $i=j-h+q$. Then
by Equations~(\ref{s61}) and (\ref{s62}), the labels on the
other edges not in $H$ incident to
$y_{2j}$ and $y_{(m_j-1)j}$,  respectively,
are
\begin{eqnarray*}
(s-2)m+2h+1-(2h+2j-1) &=& (s-2)m-2(i+h-q)+2\\
                               &=&(s-2)m-2h+2q-2i+2=d+2q-2(i-1),\\
  (s-2)m+2h+1-(2h+2j) &=& (s-2)m-2(i+h-q)+1\\
                             &=&(s-2)m-2h+2q-2i+1=d+2q-2(i-1)-1.
\end{eqnarray*}
Thus,
$$
 \omega_G(y_{2j})=\omega_H(y_{2j})+(d+2q-2(i-p-1)) , \quad
 \omega_G(y_{(m_j-1)j})=\omega_H(y_{(m_j-1)j})+(d+2q-2(i-p-1)-1).
$$
Apply Lemma~\ref{pathcycle} on $H$ with $c$ and $d$
defined above\,(according to if $s=4$ or $s\ge 6$) and with $p=0$, we get an assignment of  labels  on $E(H)$
such that
\begin{enumerate}[(i)]
 \item For any $x\in V(H)\cap X$, $\omega_H(x)=\frac{d_H(x)(c+d)}{2}$; and
 \item For any distinct $y,z\in V(H)\cap Y$, $\omega_{G}(y)\ne \omega_G(z)$,
and $\omega_{G}(y)\in [2c,2d+2q]$.
\end{enumerate}

%
%
%
%
%
%
Apply Lemma~\ref{cycle1} on $F$
with
$$
a:=(s-2)m+1 \quad \mbox{and} \quad b:=sm,
$$
we get an antimagic labeling on $F$.

If $s=4$, orient the edges in $\{x_{1i}y_{1i}, x_{1i}y_{m_ii}\,|\, 1\le i \le h\}$
from $Y$ to $X$, and
orient all the remaining edges  from $X$ to $Y$.
If $s\ge 6$,  orient the edges in $\{x_{1i}y_{1i}|\, 1\le i \le h\}$
from $Y$ to $X$, and
orient all the remaining edges  from $X$ to $Y$.
Denote the  orientation of $G$ by
$\overrightarrow{G}$.

\bigskip
{\noindent \bf Claim 3:}  The labeling of $\overrightarrow{G}$
given above  is antimagic.

\pf We first show that  the set of labels used is the set $[1,sm]$.
The labels used on edges in $F$ are exactly numbers in the
set $[(s-2)m+1, sm]$. If $s=4$,
then the set of labels used on edges in $\{x_{1i}y_{1i}, x_{1i}y_{m_ii}, x_{2i}y_{1i},
x_{m_ii}y_{m_ii}, x_{2i}y_{2i}, x_{m_ii}y_{(m_i-1)i}\,|\, 1\le i\le h\}$
is $[1,3h]\cup [2m-3h+1,2m]$, and the set of labels used on $E(H)$
is $[3h+1,2m-3h]$.
If $s\ge 6$,
then the set of labels used on edges in $\{x_{1i}y_{1i}, x_{1i}y_{m_ii}, x_{2i}y_{1i},
x_{m_ii}y_{m_ii}, x_{2i}y_{2i}, x_{m_ii}y_{(m_i-1)i}\,|\, 1\le i\le h\}$
is $[1,4h]\cup [(s-2)m-2h+1,(s-2)m]$, and the set of labels used on $E(H)$
is $[4h+1,(s-2)m-2h]$.
The union of these sets is the set  $[1,sm]$.

We then show that the oriented sums on
vertices in $\overrightarrow{G}$  are all distinct. We separate the
proof according to if $s=4$ or $s\ge 6$.

{\bf Case $s=4$:\quad }
For each $i$ with $1\le i\le h$, by~(\ref{s43})-(\ref{s46}) and the orientation of $G$,
the labels at $y_{1i}$, $y_{m_ii}$,
respectively, are
$$
(-i, h+2i-1)\quad \mbox{and} \quad (-2m+(i-1), h+2i).
$$
Thus,
$$
\omega_G(y_{1i})=h+i-1\in [h, 2h-1] \quad \mbox{and} \quad \omega_G(y_{m_ii})=-2m+h+3i-1\in[-2m+h+2,-2m+4h-1].
$$
All these $2h$ values are pairwise distinct and fall into
the interval $[-2m+h+2, 2h-1]$.

For each $i$ with $1\le i\le h-q$, by~(\ref{s41}) and (\ref{s42}),
the label at $y_{2i}$ is,
$$
(2m+1-(h+2i-1), 2m+1-(h+2i).
$$
Thus,
$$
\omega_G(y_{2i})=4m+3-2h-4i\in [4m+3-6h+4q, 4m-1-2h].
$$
All these $h-q$ values are pairwise distinct.

The sums
at vertices in $V(H)\cap Y$ are all distinct and fall into the interval
$[2c,2d+2q]=[6h+2,4m-6h+2q]\subseteq [6h+2,4m]$\,($q\le h$) by Lemma \ref{pathcycle}.
The sums at vertices in $V(F)\cap Y$ are all
distinct and fall into the interval $[4m+3, 8m-1]$
by Lemma~\ref{cycle1}. Since
these sets
$[-2m+h+2, 2h-1],  [6h+2,4m-6h+2q], [4m+3-6h+4q, 4m-1-2h]$,
and $[4m+3, 8m-1]$
are pairwise disjoint,
we see that the oriented sums on
vertices in $Y$ are all distinct.

Next we show that in $\overrightarrow{G}$,  the oriented sums on
vertices in $X$ are all distinct.
For  each $i$ with $1 \le i\le h$,
$$
\begin{cases}
     \omega_{G-E(F)}(x_{1i})=2m+1 & \text{by (\ref{s43}) and (\ref{s44}) }, \\
      \omega_{G-E(F)}(x_{2i})=\omega_{G-E(F)}(x_{m_ii})=-2m-1& \text{by (\ref{s45}) (\ref{s41}), and (\ref{s46}) (\ref{s42})},\\
      \omega_{G-E(F)}(x)=-\frac{d_H(x)(c+d)}{2}=-(c+d)=-2m-1 & \text{if $x\in V(H)\cap X$}.
\end{cases}
$$
Hence,
\begin{equation}\label{abdiffx}
|\omega_{G-E(F)}(u)-\omega_{G-E(F)}(v)|=0~ \mbox{or}~ 4m+2, \mbox{for any $ u,v\in X$}.
\end{equation}
For the graph $F$, by Lemma~\ref{cycle1}, the sums on vertices in $V(F)\cap X$  are pairwise
distinct. Since the set of labels used on $E(F)$ is $[2m+1, 4m]$ and $F$ is 2-regular,
 it follows that in $F$,  $\omega_{F}(x)\in [-8m+1,-4m-3]$ and
 any two of the sums at vertices in $V(F)\cap X$ differ an absolute  value of at most $4m-4$.
Because of $\omega_{G}(x)=\omega_{G-E(F)}(x)+\omega_{F}(x)$ for $x\in X$ and the fact in~(\ref{abdiffx}),
we conclude that the total oriented vertex sums at
vertices in $X$ are all distinct.

Finally, we show that for any $x\in X$ and $y\in Y$
the oriented vertex  sums at $x$ and $y$ are distinct in $\overrightarrow{G}$.
By the analysis above, $\omega_G(y)\in [-2m+h+2, 8m-1]$ for any $y\in Y$.
And $\omega_G(x)\in [-10m, -2m-2]$ for any $x\in X$, which follows by
the facts that $\omega_{G-E(F)}(x)=-2m-1$ or $2m+1$, $\omega_{F}(x)\in [-8m+1,-4m-3]$,
 and $\omega_G(x)=\omega_{G-E(F)}(x)+\omega_{F}(x)$.
 Thus, $\omega_G(x)\ne \omega_G(y)$.

{\bf Case $s\ge 6$}:\quad
For each $i$ with $1\le i\le h$, by~(\ref{s63})-(\ref{s66}) and the orientation of $G$, the labels at $y_{1i}$, $y_{m_ii}$,
respectively, are
$$
(-i, 2h+2i+1)\quad \mbox{and} \quad (h+i, 2h+2i).
$$
Thus,
$$
\omega_G(y_{1i})=2h+i+1\in [2h+2,3h+1] \quad \mbox{and} \quad \omega_G(y_{m_ii})=3h+3i\in[3h+3,6h].
$$
All these $2h$ values are pairwise  distinct and fall into
the interval $[2h+2, 6h]$.

For each $i$ with $1\le i\le h-q$, by~(\ref{s61}) and (\ref{s62}),
the label at $y_{2i}$ is,
$$
((s-2)m+2h+1-(h+2i-1), (s-2)m+2h+1-(h+2i).
$$
Thus,
$$
\omega_G(y_{2i})=2(s-2)m+3+2h-4i\in [ 2(s-2)m+3-2h+4q, 2(s-2)m+2h-1].
$$
All these $h-q$ values are pairwise distinct.

The sums
at vertices in $V(H)\cap Y$ are all  distinct
and fall into the interval $[2c,2d+2q]=[8h+2,2(s-2)m-4h+2q]$ by Lemma \ref{pathcycle}.
The sums at vertices in $V(F)\cap Y$ are all
distinct and fall into the interval $[2(s-2)m+3, 2sm-1]$
by Lemma~\ref{cycle1}. Since
these sets
$[2h+2, 6h], [8h+2,2(s-2)m-4h+2q]$, $[ 2(s-2)m+3-2h+4q, 2(s-2)m+2h-1]$,
and $[2(s-2)m+3, 2sm-1]$
are pairwise disjoint,
we see that the oriented vertex sums on
vertices in $Y$ are all distinct.

Next we show that in $\overrightarrow{G}$,  the oriented  sums at
vertices in $X$ are all distinct.
Assume that for each $x\in X$, $x$ appears  $\alpha_x$ times in $\{x_{1i}\,|\, 1\le i\le h\}$,
and $\beta_x$ times in  $\{x_{2i}, x_{m_ii}\,|\, 1\le i\le h\}$.
Since each of $x_{1i},x_{2i}$, and $x_{m_ii}$ has two distinct neighbors in $Y$, $x$
has degree $s-2-2\alpha_x-2\beta_x$ in $H$.
In addtion,  each appearance of $x$ in $\{x_{1i}\,|\, 1\le i\le h\}$
contributes a value of $-h$ to the oriented sum at $x$ by (\ref{s63}) and (\ref{s64}), and
each appearance of $x$ in $\{x_{2i}, x_{m_ii}\,|\, 1\le i\le h\}$ contributes a value of
$-(c+d)$ to the oriented sum at $x$ by (\ref{s65}) (\ref{s61}), and (\ref{s66}) (\ref{s62}).  By (i) preceeding Claim 3,
$\omega_{H}(x)=-\frac{d_H(x)(c+d)}{2}=-\frac{(s-2-2\alpha_x-2\beta_x)(c+d)}{2}$. Hence,
$$
\omega_{G-E(F)}(x)=-\frac{(s-2-2\alpha_x-2\beta_x)(c+d)}{2}-\alpha_xh-\beta_x(c+d)=-\frac{(s-2)(c+d)}{2}+\alpha_x(c+d-h).
$$
Thus, for any $ u,v\in X$,
\begin{eqnarray}
  |\omega_{G-E(F)}(u)-\omega_{G-E(F)}(v)| &=&|\alpha_u-\alpha_v|(c+d-h) \nonumber \\
   &=&|\alpha_u-\alpha_v|((s-2)m+h+1)=0 ~\mbox{or}~ >4m. \label{abdiffx2}
\end{eqnarray}
For the graph $F$, by Lemma~\ref{cycle1}, the sums on vertices in $V(F)\cap X$  are pairwise
distinct, any two of the sums at vertices in $V(F)\cap X$ differ an absolute  value of at most $4m-4$.
Because of $\omega_{G}(x)=\omega_{G-E(F)}(x)+\omega_{F}(x)$ for $x\in X$ and the fact in~(\ref{abdiffx2}),
we conclude that the total oriented sums at
vertices in $X$ are all distinct.

Finally, we show that for any $x\in X$ and $y\in Y$
the oriented sums at $x$ and $y$ are distinct in $\overrightarrow{G}$.
By the analysis above, for any $y\in Y$, $\omega_G(y)\in [2h+2, 2sm-1]$ is a positive integer.
For any $x\in X$,  $\omega_{G-F}(x)$ is negative and $\omega_F(x)$ is negative,
so $\omega_G(x)=\omega_{G-E(F)}(x)+\omega_{F}(x)$ is negative.
Hence  $\omega_G(x)\ne \omega_G(y)$.
\qed

The proof of Theorem~\ref{th1} is now complete.
\qed

\bibliographystyle{plain}
\bibliography{SSL-BIB}

\end{document}